# HYPOELLIPTICITY IN INFINITE DIMENSIONS AND AN APPLICATION IN INTEREST RATE THEORY[1]

By Fabrice Baudoin and Josef Teichmann

*Vienna University of Technology*

We apply methods from Malliavin calculus to prove an infinite-dimensional version of Hörmander's theorem for stochastic evolution equations in the spirit of Da Prato–Zabczyk. This result is used to show that HJM-equations from interest rate theory, which satisfy the Hörmander condition, have the conceptually undesirable feature that any selection of yields admits a density as multi-dimensional random variable.

**1. Introduction.** Given a separable Hilbert space $H$ and the generator $A$ of a strongly continuous group (sic!), we aim to prove a Hörmander theorem for stochastic evolution equations of the Da Prato–Zabczyk type (see [4] for all details)

$$dr_t = (Ar_t + \alpha(r_t)) \, dt + \sum_{i=1}^{d} \sigma_i(r_t) \, dB_t^i,$$

$(E_0)$

$$r_0 \in H,$$

under the assumption that iterative Lie brackets of the Stratonovich drift and the volatility vector fields span the Hilbert space. We therefore apply methods from Malliavin calculus, which have already been used to solve similar questions in filtering theory (see, e.g., [12]) in stochastic differential geometry (see, e.g., [1] and [2]) or in stochastic analysis (see, e.g., [7]).

A particular example, which received some attention recently (see, e.g., [3] and [6]) is the Heath–Jarrow–Morton equation of interest rate theory (in the sequel abbreviated by HJM),

$$dr_t = \left( \frac{d}{dx} r_t + \alpha_{\mathrm{HJM}}(r_t) \right) dt + \sum_{i=1}^{d} \sigma_i(r_t) \, dB_t,$$

Received March 2004; revised September 2004.
[1] Supported by the Research Training Network HPRN-CT-2002-00281.

where $H$ is a Hilbert space of real-valued functions on the real line. The HJM drift term is given by

$$\alpha_{\mathrm{HJM}}(r)(x) := \sum_{i=1}^{d} \sigma^i(r)(x) \int_0^x \sigma^i(r)(y) \, dy$$

for $x \geq 0$ and $r \in H$. In order to apply Theorem 1 to the HJM-equation, we introduce the relevant setting in Section 3.

The HJM-equation describes the time-evolution of forward rates (which contain the full information of a considered bond market) in the martingale measure. It is of particular importance in applications to identify relevant, economically reasonable factors in this evolution. More precisely, how do you find a Markov process with values in some finite-dimensional state space (the space of economically reasonable factors), such that the whole evolution becomes a deterministic function of this Markov process? Conditions in order to guarantee this behavior have been described in [3] and [6]. Economically reasonable factors are the forward rate itself at some time to maturity $x \geq 0$, or averages drawn from it, so-called Yields. If the time-evolution of interest rates cannot be described by finitely many stochastic factors, we can imagine the following generic behavior, which we formulate in a criterion.

CRITERION 1.   We denote by $(r_t(x))_{t \geq 0}$ a forward rate evolution in the Musiela parametrization, that is, a mild solution of the HJM equation. For $x > 0$, the associated Yield is denoted by

$$Y_t(x) := \frac{1}{x} \int_0^x r_t(y) \, dy,$$

we define $Y_t(0) = R_t = r_t(0)$, the short rate process. The evolution is called *generic* if for each selection of times to maturity $0 \leq x_1 < \cdots < x_n$, the $\mathbb{R}^n$-valued process $(Y_t(x_1), \ldots, Y_t(x_n))$ admits a density with respect to the Lebesgue measure.

REMARK 1.   In financial mathematics generic evolutions do not seem reasonable, since—loosely speaking—the support of the random variable $r_t$, for some $t > 0$, becomes too big in the Hilbert space of forward rate curves. In other words, any "shape" of forward rate curves, which we assume from the beginning to model the market phenomena, is destroyed with positive probability. Hence, the very restrictive phenomenon of finite-dimensional realizations for the HJM-equation also appears as the only structure where "shape" is not destroyed immediately. Hence, generic evolutions behave essentially different from affine, finite-dimensional realizations, where we can always find tenors $x_1 < \cdots < x_n$, such that the Yield process does not admit a density.

By [6], the existence of finite-dimensional realizations is—among technical assumptions—equivalent to the fact that the stochastic evolution admits locally invariant submanifolds (with boundary). This is equivalent to the fact that a certain Lie algebra of vector fields $\mathcal{D}_{LA}$ is evaluated to a finite-dimensional subspace of the Hilbert space $H$ at "some" points $r \in H$, more precisely, there is a natural number $M \geq 1$, such that

$$\dim_{\mathbb{R}} \mathcal{D}_{LA}(r) \leq M < \infty$$

in a $\operatorname{dom}(A^\infty)$-neighborhood.

In Section 2 we prove the Hörmander-type result for evolution equations where the drift contains a group generator. We then show in Section 3 that for generic volatility structures at a point $r_0 \in H$, the HJM-equation leads to a generic evolution for the initial value $r_0$.

Conceptually, a generic evolution is not desirable in interest rate theory, since we expect to exhaust all information by a finite number of Yields. Hence, the result Theorem 2 can be interpreted as an additional argument for finite-dimensional realizations. Notice also that this result is invariant under the important equivalent changes of measure: if we obtain a generic evolution with respect to one fixed measure, then also with respect to all equivalent measures.

**2. Malliavin calculus in Hilbert spaces.** In order to set up the methodological background, we refer, on the one hand, to the finite-dimensional literature in Malliavin calculus, such as [11]. On the other hand, we refer to [6] for the analytical framework, in particular, for questions of differentiability of functions on infinite-dimensional spaces and for the notion of derivatives of vector fields $V : U \subset G \to G$, when $G$ is some Fréchet space.

We shall mainly work on Hilbert spaces: then the derivative $DV : U \to L(H)$ is a linear operator to the *Banach space* of bounded linear operators, where we can speak about usual properties as differentiability, boundedness, and so on.

We consider evolution equations of the type

$$(E_0) \quad \begin{aligned} dr_t &= (Ar_t + \alpha(r_t))\,dt + \sum_{i=1}^{d} \sigma_i(r_t)\,dB_t^i, \\ r_0 &\in H, \end{aligned}$$

where $A : \operatorname{dom}(A) \subset H \to H$ is the generator of a strongly continuous group $(T_t)_{t \geq 0}$ on a separable Hilbert space $H$. We apply, furthermore, the following notation:

$$\operatorname{dom}(A^k) := \{h \in H \mid h \in \operatorname{dom}(A^{k-1}) \text{ and } A^{k-1}h \in \operatorname{dom}(A)\},$$

$$\|h\|^2_{\mathrm{dom}(A^k)} := \sum_{i=0}^{k} \|A^i h\|^2,$$

$$\mathrm{dom}(A^\infty) = \bigcap_{k \geq 0} \mathrm{dom}(A^k).$$

The maps $\alpha, \sigma_1, \ldots, \sigma_d : H \to \mathrm{dom}(A^\infty)$ are smooth vector fields with the property that $\alpha, \sigma_1, \ldots, \sigma_d : \mathrm{dom}(A^k) \to \mathrm{dom}(A^k)$ are $C^\infty$-bounded. As usual (see, e.g., [11]), a vector field $V$ is called $C^\infty$-*bounded* if each higher derivative $D^l V : \mathrm{dom}(A^k) \to L^l(\mathrm{dom}(A^k))$ is a bounded function for $l \geq 1$. In this case $V$ grows at most linearly on $\mathrm{dom}(A^k)$.

Notice that due to the regularity assumptions, we can interpret the equation ($E_0$) also on the Hilbert space $\mathrm{dom}(A^k)$, with the same regularity conditions on $C^\infty$-boundedness,

$$dr_t = (Ar_t + \alpha(r_t)) \, dt + \sum_{i=1}^{d} \sigma_i(r_t) \, dB_t^i,$$

($E_k$)

$$r_0 \in \mathrm{dom}(A^k).$$

A global, mild, continuous solution of equation ($E_k$) with initial value $r_0 \in \mathrm{dom}(A^k)$ is an adaped stochastic process with continuous paths $(r_t)_{t \geq 0}$ such that

$$r_t = T_t r_0 + \int_0^t T_{t-s} \alpha(r_s) \, ds + \sum_{i=1}^{d} \int_0^t T_{t-s} \sigma_i(r_s) \, dB_s^i$$

for $t \geq 0$, where $T$ is the group generated by $A$. Clearly, every strong, continuous solution is a mild, continuous solution by variation of constants (see [4]). We shall often use the vector field $\mu$, referred to as Stratonovich drift,

$$\mu(r) := Ar + \alpha(r) - \frac{1}{2} \sum_{i=1}^{d} D\sigma_i(r) \cdot \sigma_i(r)$$

for $r \in \mathrm{dom}(A)$. Notice that the Stratonovich drift is only well defined on a dense subspace $\mathrm{dom}(A^{k+1})$ of $\mathrm{dom}(A^k)$ for $k \geq 0$, if we want $\mu$ to take values in $\mathrm{dom}(A^k)$. Furthermore, $\mu$ is not even continuous. We, nevertheless, have the following regularity result:

PROPOSITION 1. *Given equation ($E_k$), for every $r_0 \in \mathrm{dom}(A^k)$, there is a unique, global mild solution with continuous paths denoted by $(r_t)_{t \geq 0}$. The natural injections $\mathrm{dom}(A^k) \to \mathrm{dom}(A^{k+1})$ leave solutions invariant, that is, a solution of equation ($E_k$) with initial value in $\mathrm{dom}(A^{k+1})$ is a also a solution of the equation with index $k + 1$. More precisely, a mild solution with initial value in $\mathrm{dom}(A^{k+1})$ is also a mild solution of the equation with index $k + 1$, and, hence, a strong solution of equation ($E_k$).*

A mild solution of equation $(E_k)$ with initial value $r_0 \in \mathrm{dom}(A^{k+1})$ is a strong solution of equation $(E_k)$, hence, the solution process is a semimartingale and the Stratonovich decomposition makes sense,

$$dr_t = \mu(r_t)\,dt + \sum_{i=1}^{d} \sigma_i(r_t) \circ dB_t^i.$$

If we assume that $r_0 \in \mathrm{dom}(A^\infty)$, then we can construct a solution process $(r_t)_{t \geq 0}$ with continuous trajectories in $\mathrm{dom}(A^\infty)$, since the Picard–Lindelöf approximation procedure converges in every Hilbert space $\mathrm{dom}(A^k)$, and the topology of $\mathrm{dom}(A^\infty)$ is the projective limit of the ones on $\mathrm{dom}(A^k)$.

For equations of the above type, the following regularity assertions hold true for the first variation process.

PROPOSITION 2. *The first variation equations, with respect to $(E_k)$ for $k \geq 0$, are well defined on $\mathrm{dom}(A^k)$*

$$dJ_{s \to t}(r_0) \cdot h = (A(J_{s \to t}(r_0) \cdot h) + D\alpha(r_t) \cdot J_{s \to t}(r_0) \cdot h)\,dt$$
$$+ \sum_{i=1}^{d} D\sigma_i(r_t) \cdot (J_{s \to t}(r_0) \cdot h)\,dB_t^i,$$
$$J_{s \to s}(r_0) \cdot h = h,$$

*for $h \in \mathrm{dom}(A^k)$, $r_0 \in \mathrm{dom}(A^k)$ and $k \geq 0$, $t \geq s$. The Stratonovich decomposition on $\mathrm{dom}(A^k)$*

(2.1)
$$dJ_{s \to t}(r_0) \cdot h = D\mu(r_t) \cdot (J_{s \to t}(r_0) \cdot h)\,dt$$
$$+ \sum_{i=1}^{d} D\sigma_i(r_t) \cdot (J_{s \to t}(r_0) \cdot h) \circ dB_t^i$$

*is only well defined for $h, r_0 \in \mathrm{dom}(A^{k+1})$, since we need to integrate semimartingales. The Itô equation has unique global mild solutions and $J_{s \to t}(r_0)$ defines a continuous linear operator on $\mathrm{dom}(A^k)$, which is invertible if $r_0 \in \mathrm{dom}(A^{k+1})$, $k \geq 0$. The adjoint of the inverse $(J_{s \to t}(r_0)^{-1})^*$ admits the Stratonovich decomposition*

(2.2)
$$d(J_{s \to t}(r_0)^{-1})^* \cdot h$$
$$= -D\mu(r_t)^* \cdot ((J_{s \to t}(r_0)^{-1})^* \cdot h)\,dt$$
$$- \sum_{i=1}^{d} D\sigma_i(r_t)^* \cdot (J_{s \to t}(r_0)^{-1})^* \cdot h \circ dB_t^i$$

*for $h, r_0 \in \mathrm{dom}(A^{k+1})$ and $k \geq 0$, $t \geq s \geq 0$. We have, furthermore,*

$$J_{s \to t}(r_0) = J_{0 \to t}(r_0) J_{0 \to s}(r_0)^{-1}$$

$\mathbb{P}$*-almost surely for $t \geq s \geq 0$.*

REMARK 2. We define a Hilbert space $\mathcal{H}_k([0,T])$ of progressively measurable processes $(r_s)_{0 \leq s \leq T}$ such that

$$\mathbb{E}\left(\int_0^T \|r_s\|_{\mathrm{dom}(A^k)}^2\right) < \infty.$$

Solutions of equations $(E_k)$ can be viewed as mappings $r_0 \mapsto (r_t)_{0 \leq t \leq T}$. Then $J_{0 \rightarrow T}(r_0) \cdot h$ is the derivative of this map in the respective locally convex structures.

REMARK 3. For the proof of Proposition 2, we need the property that $A$ and $-A$ generate a semi-group, which is equivalent to the assertion that $A$ generates a strongly continuous group (see [4] for further references). Under this assumption, we can solve the equations for $J_{0 \rightarrow s}(r_0) \cdot h$ in the Hilbert spaces $\mathrm{dom}(A^k)$ and obtain invertibility as asserted. If $A$ does not generate a strongly continuous group, the first variations will not be invertible in general.

PROOF OF PROPOSITION 2. Under our assumptions, the regularity assertions are clear, also the calculation of the first variations (see [4] for all details). The only point left is that we are allowed to pass to the Stratonovich decomposition, which is correct, since the assertions of Proposition 1 apply and since we integrate semi-martingales by Itô's formula on Hilbert spaces (see [4]). Fix now $r_0 \in \mathrm{dom}(A^{k+1})$ and $h \in \mathrm{dom}(A^{k+1})$, then invertibility follows from the fact that the semi-martingale

$$(\langle J_{s \rightarrow t}(r_0) \cdot h_1, (J_{s \rightarrow t}(r_0)^{-1})^* \cdot h_2 \rangle_{\mathrm{dom}(A^k)})_{t \geq s \geq 0}$$

is constant by the respective Stratonovich decomposition, which leads to

$$\langle J_{s \rightarrow t}(r_0)^{-1} \cdot J_{s \rightarrow t}(r_0) \cdot h_1, h_2 \rangle_{\mathrm{dom}(A^k)} = \langle h_1, h_2 \rangle_{\mathrm{dom}(A^k)}$$

for $h_1, h_2 \in \mathrm{dom}(A^{k+1})$. From this, we obtain left invertibility by continuity.

To prove that the left inverse also is a right inverse, we shall apply the following reasoning. Given an ortho-normal basis $(g_i)_{i \geq 1}$ of $\mathrm{dom}(A^k)$ which lies in $\mathrm{dom}(A^{k+1})$, we can easily compute the semi-martingale decomposition of

$$\sum_{i=1}^N \langle J_{s \rightarrow t}(r_0)^{-1} \cdot h_1, g_i \rangle_{\mathrm{dom}(A^k)} \langle g_i, J_{s \rightarrow t}(r_0)^* \cdot h_2 \rangle_{\mathrm{dom}(A^k)}$$

$$= \sum_{i=1}^N \langle h_1, (J_{s \rightarrow t}(r_0)^{-1})^* \cdot g_i \rangle_{\mathrm{dom}(A^k)} \langle J_{s \rightarrow t}(r_0) \cdot g_i, h_2 \rangle_{\mathrm{dom}(A^k)},$$

for $h_1, h_2 \in \mathrm{dom}(A^{k+1})$ and $N \geq 1$. Now we apply the Stratonovich decomposition: by adjoining, we can free the $g_i$'s and pass to the limit, which yields vanishing finite variation and martingale part. Hence,

$$\langle J_{s\to t}(r_0)J_{s\to t}(r_0)^{-1}\cdot h_1, h_2\rangle_{\mathrm{dom}(A^k)}$$
$$= \lim_{N\to\infty}\sum_{i=1}^{N}\langle J_{s\to t}(r_0)^{-1}\cdot h_1, g_i\rangle_{\mathrm{dom}(A^k)}\langle g_i, J_{s\to t}(r_0)^*\cdot h_2\rangle_{\mathrm{dom}(A^k)}$$
$$= \langle h_1, h_2\rangle_{\mathrm{dom}(A^k)},$$

which is the equation for the right inverse.

Finally, the process $(J_{0\to t}(r_0)J_{0\to s}(r_0)^{-1})_{t\geq s}$ satisfies the correct differential equation and we obtain by uniqueness the desired assertion on the decomposition of the first variation process $J_{s\to t}(r_0)$. □

Crucial for the further analysis is the notion of the Lie bracket of two vector fields $V_1, V_2 : \mathrm{dom}(A^\infty)\to\mathrm{dom}(A^\infty)$ (see [6] for the analytical framework). We need to leave the category of Hilbert spaces, since the vector field $\mu$ is only well defined on $\mathrm{dom}(A^\infty)$ as a smooth vector field. We define

$$[V_1, V_2](r) := DV_1(r)\cdot V_2(r) - DV_2(r)\cdot V_1(r)$$

for $r\in\mathrm{dom}(A^\infty)$.

Fix $r_0\in\mathrm{dom}(A^\infty)$. We then define the distribution $\mathcal{D}(r_0)\subset H$, which is generated by $\sigma_1(r_0),\ldots,\sigma_d(r_0)$ and their iterative Lie brackets with the vector fields $\mu, \sigma_1,\ldots,\sigma_d$, evaluated at the point $r_0$. Notice that a priori the direction $\mu(r_0)$ does not appear in the definition of $\mathcal{D}(r_0)$,

$$\mathcal{D}(r_0) = \langle\sigma_1(r_0),\ldots,\sigma_d(r_0),[\sigma_i,\sigma_j](r_0),\ldots,[\mu,\sigma_i](r_0),\ldots\rangle.$$

As in the finite-dimensional case, following the original idea of Malliavin [9], the main theorem is proved by calculation of the (reduced) covariance matrix (see also [11] for a more recent presentation). We need an additional lemma on Lie brackets of the type $[\mu,\sigma_i]$ for $i=1,\ldots,d$.

LEMMA 1. *Given a vector field* $V : \mathrm{dom}(A^k)\to\mathrm{dom}(A^\infty)$, *then there is a smooth extension of the Lie bracket* $[\mu, V] : \mathrm{dom}(A^{k+1})\to\mathrm{dom}(A^\infty)$.

PROOF. A vector field $V : \mathrm{dom}(A^k)\to\mathrm{dom}(A^\infty)$ is well defined and smooth on $\mathrm{dom}(A^\infty)\subset\mathrm{dom}(A^k)$. There we define the Lie bracket with $\mu$ and obtain a well-defined vector field $[\mu, V] : \mathrm{dom}(A^\infty)\to\mathrm{dom}(A^\infty)$. Take $\mu(r) = Ar + \beta(r)$, where $\beta : \mathrm{dom}(A^k)\to\mathrm{dom}(A^\infty)$. Then

$$[\mu, V](r) := AV(r) + D\beta(r)\cdot V(r) - DV(r)\cdot Ar - DV(r)\cdot\beta(r).$$

Since $DV(r) : \mathrm{dom}(A^k)\to\mathrm{dom}(A^\infty)$, we obtain a smooth extension on $\mathrm{dom}(A^{k+1})$. For details and further references on the analysis, see [6]. □

THEOREM 1. *Fix $r_0 \in \mathrm{dom}(A^\infty)$ and assume that $\mathcal{D}(r_0)$ is dense in $H$. Given $k$ linearly independent functionals $\ell := (l_1, \ldots, l_k) : H \to \mathbb{R}$, the law of the process $(\ell \circ r_t)_{t \geq 0}$ admits a density with respect to Lebesgue measure on $\mathbb{R}^k$ for $t > 0$.*

PROOF. Take $t > 0$. We have to form the Malliavin covariance matrix $\gamma_t$, which is done by well-known formulas on the first variation (see [11]). The covariance matrix can be decomposed into

$$\gamma_t = (\ell \circ J_{0 \to t}(r_0)) C_t (\ell \circ J_{0 \to t}(r_0))^\mathsf{T},$$

where the random, symmetric Hilbert–Schmidt-operator $C_t$, the reduced covariance operator, is defined via

$$\langle y, C_t y \rangle = \sum_{p=1}^{d} \int_0^t \langle y, J_{0 \to s}(r_0)^{-1} \cdot \sigma_p(r_s) \rangle^2 \, ds.$$

We first show that $C_t$ is a positive operator. We denote the kernel of $C_t$ by $K_t \subset H$ and get a decreasing sequence of closed random subspaces of $H$. $V = \bigcup_{t > 0} K_t$ is a deterministic subspace by the Blumenthal zero–one law, that is, there exists a null set $N$ such that $V$ is deterministic on $N^c$. We shall do the following calculus on $N^c$.

We fix $y \in V$, then we consider the stopping time

$$\theta := \inf\{s, q_s > 0\},$$

with respect to the continuous semi-martingale

$$q_s = \sum_{p=1}^{d} \langle y, J_{0 \to s}(r_0)^{-1} \cdot \sigma_p(r_s) \rangle^2.$$

Then $\theta > 0$ almost surely and $q_{s \wedge \theta} = 0$ for $s \geq 0$.

Now, a continuous $L^2$-semi-martingale with values in $\mathbb{R}$

$$s_s - s_0 = \sum_{k=1}^{d} \int_0^s \alpha_k(u) \, dB_u^k + \int_0^s \beta(u) \, du$$

for $s \geq 0$, which vanishes up to the stopping time $\theta$, satisfies—due to the Doob–Meyer decomposition—

$$\alpha_k(s \wedge \theta) = 0$$

for $k = 1, \ldots, d$.

We shall apply this consideration for the continuous semi-martingales $m_s := \langle y, J_{0 \to s}(r_0)^{-1} \cdot \sigma_p(r_s) \rangle$ on $[0, t]$ for $p = 1, \ldots, d$. Therefore, we need to calculate the Doob–Meyer decomposition of $(m_s)_{0 \leq s \leq t}$. This can be done simply by applying equation (2.2) for the adjoint of $J_{0 \to s}(r_0)^{-1}$,

$$
\begin{aligned}
dm_s &= \langle d(J_{0 \to s}(r_0)^{-1})^* \cdot y, \sigma_p(r_s) \rangle \\
&= -\langle D\mu(r_s)^* \cdot (J_{0 \to s}(r_0)^{-1})^* \cdot y, \sigma_p(r_s) \rangle \, ds \\
&\quad - \sum_{i=1}^{d} \langle D\sigma_i(r_s)^* \cdot (J_{0 \to s}(r_0)^{-1})^* \cdot y, \sigma_p(r_s) \rangle \circ dB_s^i \\
&\quad + \langle (J_{0 \to s}(r_0)^{-1})^* \cdot y, D\sigma_p(r_s) \cdot \mu(r_s) \rangle \, ds \\
&\quad + \sum_{i=1}^{d} \langle (J_{0 \to s}(r_0)^{-1})^* \cdot y, D\sigma_p(r_s) \cdot \sigma_i(r_s) \rangle \circ dB_s^i \\
&= \langle (J_{0 \to s}(r_0)^{-1})^* \cdot y, [\sigma_p, \mu](r_s) \rangle \, ds \\
&\quad + \sum_{i=1}^{d} \langle (J_{0 \to s}(r_0)^{-1})^* \cdot y, [\sigma_p, \sigma_i](r_s) \rangle \circ dB_s^i.
\end{aligned}
$$

From the Doob–Meyer decomposition this leads to

$$
\begin{aligned}
\langle y, J_{0 \to s}(r_0)^{-1} \cdot [\sigma_p, \sigma_i](r_s) \rangle &= 0, \\
\langle y, J_{0 \to s}(r_0)^{-1} \cdot [\sigma_p, \mu](r_s) \rangle &= 0
\end{aligned}
$$

for $i = 1, \ldots, d$, $p = 1, \ldots, d$ and $0 \leq s \leq \theta$. Notice that all the appearing Lie brackets have a smooth extension to some $\mathrm{dom}(A^k)$ for $k \geq 0$ due to Lemma 1, where we can repeat the argument recursively.

Consequently, the above equation leads by iterative application to

$$
\langle y, J_{0 \to s}(r_0)^{-1} \cdot \mathcal{D}(r_s) \rangle = 0
$$

for $s \leq \theta$. Evaluation at $s = 0$ yields $y = 0$, since $\mathcal{D}(r_0)$ is dense in $H$. Hence, $C_t$ is a positive definite operator. Therefore, we obtain that there is a null set $N$, such that on $N^c$ the matrix $C_t$ has an empty kernel. Hence, the law is absolutely continuous with respect to Lebesgue measure, since $J_{0 \to t}(r_0)$ is invertible and, therefore, $\gamma_t$ has an empty kernel (Theorem 2.1.2 in [11], page 86). □

EXAMPLE 1. For instance, if we consider the equation

$$
dr_t = Ar_t \, dt + \sum_{i=1}^{d} h_i \, dB_t^i,
$$

where $h_1, \ldots, h_d \in \mathrm{dom}(A^\infty)$, then it is easily seen that $\mathcal{D}(r_0)$ is dense in $H$ as soon as the linear span of the orbit $(A^n h_i)_{n \geq 0, 1 \leq i \leq d}$ is dense in $H$, for all $r_0 \in \mathrm{dom}(A^\infty)$. As an example, we can consider $H = L^2(\mathbb{R}, \lambda)$, where $\lambda$ denotes the Lebesgue measure on $\mathbb{R}$, $d = 1$, $A = \frac{d}{dx}$ and $h_1 = e^{-x^2/2}$. This result is well known and can be obtained by simpler methods.

EXAMPLE 2. For non-Gaussian random variables the assertion of Theorem 1 is already nontrivial. Let $\sigma(r) = \phi(r)h$ be a smooth vector field, $\phi \colon H \to \mathbb{R}$ a $C^\infty$-bounded function, fix $r_0 \in \mathrm{dom}(A^\infty)$ such that $\phi(r_0) \neq 0$. Then we can calculate conditions such that the Lie algebra at $r_0$ is dense in $H$,

$$\mu(r) = Ar - \tfrac{1}{2}(\phi(r)D\phi(r) \cdot h)h$$
$$= Ar + \psi(r)h,$$
$$D\mu(r) \cdot g = Ag + (D\psi(r) \cdot g)h,$$
$$[\mu, \sigma](r) = \phi(r)Ah + \phi(r)(D\psi(r) \cdot h)h - (D\phi(r) \cdot Ar)h - \psi(r)(D\phi(r) \cdot h)h$$
$$= \phi(r)Ah + \psi_2(r)h,$$

hence, the span of $A^n h$ lies in $\mathcal{D}(r_0)$ (division by $\phi$ around $r_0$ is performed). Consequently, a necessary condition for $\mathcal{D}(r_0)$ to be dense in $H$ is that the linear span of the orbit $(A^n h)_{n \geq 0}$ is dense in $H$.

**3. Applications to interest rate theory.** We shall describe a framework for the HJM-equation, where Theorem 1 applies. This framework is narrower than the setting given in [5], but it enables us to conclude without worries the desired result.

- $H$ is a separable Hilbert space of continuous functions on the whole real line containing the constant functions (constant term structures). The point evaluations are continuous with respect to the topology of the Hilbert space. Furthermore, we assume that the long rate is well defined and a continuous linear functional $l_\infty(r) := r(\infty)$ for $r \in H$.
- The shift semigroup $(T_t r)(x) = r(t + x)$ is a strongly continuous group on $H$ with generator $\frac{d}{dx}$.
- The map $h \mapsto \mathcal{S}(h)$ with $\mathcal{S}(h)(x) := h(x)\int_0^x h(y)\,dy$ for $x \geq 0$ (if $x < 0$ this relation need not hold true) satisfies

$$\|\mathcal{S}(h)\| \leq K\|h\|^2$$

for all $h \in H$ with $\mathcal{S}(h) \in H$ for some real constant $K$. There is a closed subspace $H_0 \subset H$ such that $\mathcal{S}(h) \in H$ if and only if $h \in H_0$.

EXAMPLE 3. The first example and seminal treating of consistency problems in interest rate theory is outlined in [3]. Here the Hilbert space $H$ is a space of entire functions, where all the requirements above are fulfilled. In particular, the shift group is generated by a bounded operator $\frac{d}{dx}$ on this Hilbert space.

Example 4. Hilbert spaces of the above type can be constructed by methods similar to [5], pages 75–81, and can be chosen of the type (for the notation, see [5])

$$H_w := \left\{ h \in H^1_{\mathrm{loc}}(\mathbb{R}) \mid \|h\|^2_w := \int_{-\infty}^{\infty} |h'(x)|^2 \, w(x) \, dx + |h(0)|^2 < \infty \right\}.$$

Notice that, in contrary to [5], we need the forward rate curves to be defined on the whole real line. The forward curve on the negative real line has no direct financial interpretation.

We need a further requirement for the volatility vector fields in order that the function $\mathcal{S}$ is well defined: we define $\mathrm{dom}((\frac{d}{dx})_0) = \mathrm{dom}(\frac{d}{dx}) \cap H_0$ and similar for all further powers.

- The vector fields are smooth maps $\sigma^i \colon H \to \mathrm{dom}((\frac{d}{dx})_0^\infty)$ for $i = 1, \ldots, d$.
- The restriction $\sigma^i \colon U \cap \mathrm{dom}((\frac{d}{dx})^k) \to \mathrm{dom}((\frac{d}{dx})^k)$ are $C^\infty$-bounded for $i = 1, \ldots, d$ and $k \geq 0$.
- The HJM drift term is defined to be $\sum_{i=1}^d \mathcal{S}(\sigma^i)$

$$l_\infty \left( \sum_{i=1}^d \mathcal{S}(\sigma^i) \right) = 0,$$

where $l_\infty$ denotes the linear functional mapping a term structure to its long rate $r(\infty)$. By [8], the long rate of an arbitrage free term structure model is an increasing process, hence, this condition means that we assume it to be constant.

Under these conditions, we can prove the following lemma, which guarantees that the hypoellipticity result can be applied.

Lemma 2. *Let the above conditions be in force. Then the Hilbert space $H_0 := \ker l_\infty$ of term structures vanishing at $\infty$ is an invariant subspace of the HJM equation, furthermore,*

$$l_\infty(r_t) = r^*(\infty)$$

*is deterministic for $t \geq 0$ for any solution $(r_t)_{t \geq 0}$ with initial value $r^*$ of the HJM equation.*

Proof. We take a mild solution of the HJM equation with initial value $r^*$,

$$r_t = T_t r^* + \int_0^t T_{t-s} \alpha_{\mathrm{HJM}}(r_s) \, ds + \sum_{i=1}^d \int_0^t T_{t-s} \sigma^i(r_s) \, dB_s,$$

and apply the linear functionals $l_\infty$ to this equation. By continuity, we obtain

$$l_\infty(r_t) = l_\infty(T_t r^*) + \int_0^t l_\infty(T_{t-s}\alpha_{\mathrm{HJM}}(r_s))\,ds + \sum_{i=1}^d \int_0^t l_\infty(T_{t-s}\sigma^i(r_s))\,dB_s^i$$

$$= l_\infty(T_t r^*),$$

since $l_\infty(T_{t-s}\alpha_{\mathrm{HJM}}(r)) = 0$ and $l_\infty(T_{t-s}\sigma^i(r)) = 0$ for $r \in U \subset H$ by our assumptions. $\square$

With respect to this subspace of codimension 1, we can suppose that the condition

$$\mathcal{D}(r_0) \text{ is dense in } H_0$$

holds true, since the only deterministic direction of the time evolution, namely, $l_\infty$, is excluded.

THEOREM 2. *Take the above conditions and assume that, for some $r_0 \in \mathrm{dom}((\frac{d}{dx})^\infty)$, the condition*

$$\mathcal{D}(r_0) \text{ is dense in } H_0$$

*holds true, then for linearly independent linear functionals $l_1, \ldots, l_n : H_0 \to \mathbb{R}$, the random variable $(l_1(r_t), \ldots, l_n(r_t))$ has a density.*

PROOF. We have to restrict the reasoning to $H_0$. Take $r_0 \in H$ and define $r^* = r_0 - r_0(\infty) \in H_0$ [subtracting the constant term structure at level $r_0(\infty)$]. With the new vector fields,

$$\sigma^i(r) := \sigma^i(r + r_0(\infty))$$

for $r \in H_0$ and $i = 1, \ldots, d$. The solution of the equation associated to these vector fields at initial value $r^*$ is given through $(r_t - r_0(\infty))_{t \geq 0}$, where $(r_t)_{t \geq 0}$ denotes the solution of the original equation with initial value $r_0$. Since the Lie algebraic condition does not change under translation, we can conclude by Theorem 1 that for the given linearly independent $l_1, \ldots, l_n$, the random variable $(l_1(r_t), \ldots, l_n(r_t))$ has a density with respect to Lebesgue measure for $t > 0$. $\square$

COROLLARY 1. *Assume that*

$$\mathcal{D}(r_0) \text{ is dense in } H_0,$$

*then the evolution of the term structure of interest rates is generic.*

PROOF. For $x_1 < \cdots < x_n$, the linear functionals $Y_i(r) := \int_0^{x_i} r(y)\,dy$ for $i = 1, \ldots, n$ are linearly independent as linear functionals on the subspace $H_0$. $\square$

# REFERENCES

[1] AIRAULT, H. (1989). Projection of the infinitesimal generator of a diffusion. *J. Funct. Anal.* **85** 353–391. MR1012210

[2] BISMUT, J. M. (1981). Martingales, the Malliavin calculus, and hypoellipticity under general Hörmander conditions. *Z. Wahrsch. Verw. Gebiete* **56** 469–505. MR621660

[3] BJÖRK, T. and SVENSSON, L. (2001). On the existence of finite-dimensional realizations for nonlinear forward rate models. *Math. Finance* **11** 205–243. MR1822777

[4] DA PRATO, G. and ZABCZYK, J. (1992). *Stochastic Equations in Infinite Dimensions.* Cambridge Univ. Press. MR1207136

[5] FILIPOVIĆ, D. (2001). *Consistency Problems for Heath–Jarrow–Morton Interest Rate Models.* Springer, Berlin. MR1828523

[6] FILIPOVIĆ, D. and TEICHMANN, J. (2003). Existence of invariant manifolds for stochastic equations in infinite dimension. *J. Funct. Anal.* **197** 398–432. MR1960419

[7] HAIRER, M., MATTINGLY, J. and PARDOUX, É. (2004). Malliavin calculus for highly degenerate 2D stochastic Navier–Stokes equations. *C. R. Math. Acad. Sci. Paris* **339** 793–796. MR2110383

[8] HUBALEK, F., KLEIN, I. and TEICHMANN, J. (2002). A general proof of the Dybvig–Ingersoll–Ross theorem: Long forward rates can never fall. *Math. Finance* **12** 447–451. MR1926241

[9] MALLIAVIN, P. (1978). Stochastic calculus of variations and hypoelliptic operators. In *Proc. Inter. Symp. on Stoch. Diff. Equations, Kyoto 1976* 195–263. Wiley, New York. MR536013

[10] MALLIAVIN, P. (1997). *Stochastic Analysis.* Springer, Berlin. MR1450093

[11] NUALART, D. (1995). *The Malliavin Calculus and Related Topics.* Springer, Berlin. MR1344217

[12] OCONE, D. (1988). Stochastic calculus of variations for stochastic partial differential equations. *J. Funct. Anal.* **79** 288–331. MR953905

INSTITUTE OF MATHEMATICAL METHODS IN ECONOMICS
RESEARCH UNIT E105 FINANCIAL AND ACTUARIAL MATHEMATICS
VIENA UNIVERSITY OF TECHNOLOGY
WIEDNER HAUPTSTRASSE 8-10/105
A-1040 VIENA
AUSTRIA
E-MAIL: baudoin@fam.tuwien.ac.at
E-MAIL: jteichma@fam.tuwien.ac.at